% eb.tex  28 August 2003
\magnification 1200

\def\ver{{\it Version 1.0}}
\def\vdate{ {\it 1 September 2003}}
%\headline{\hfil \ver\hfil \vdate\hfil}

%% Macros.tex  Updated 14.8.02

   %%SPACINGS  
\def\bigb{\bigbreak\noindent}         
\def\med{\medbreak\noindent}

\def\ms{\medskip}

\def\small{\smallskip\noindent}
\def\nl{\hfil\break}

  %%TODAY'S DATE
\def\today{\noindent\number\day
\space\ifcase\month\or
  January\or February\or March\or April\or May\or June\or
  July\or August\or September\or October\or November\or December\fi
  \space\number\year}

  %%FORMATTING
\def\cl {\centerline}

\def\section#1#2{{\bigbreak\bigskip \centerline{\bf #1. #2}\bigskip}}
\def \nr { \smallskip \noindent }

  %%MATHS
\input amssym.def
\input amssym.tex

  %%Open

\def\bP {{\Bbb P}} \def\bE {{\Bbb E}} 
\def\bR {{\Bbb R}}  \def\bZ {{\Bbb Z}}

  %%Script
\def\sA {{\cal A}}  
\def\sD {{\cal D}} \def\sE {{\cal E}}

  %% Overline
\def\ol{\overline}           
\def\oB {{\overline B}}

%\def\ob{\overline B}

  %%Greek

\def\lam {\lambda}  \def\Gam{\Gamma}

  \def\eps{\varepsilon}

\def\th{\theta}

  %%Maths
\def\ra {\rightarrow}

\def\pd {\partial}
\def\q{\quad} 
\def\dint{\int\kern-.6em\int}

  %%Mathops

\def\Osc{\mathop{{\rm Osc}}}

%\def\liminf{\mathop{\underline{\rm lim}}}
%\def\limsup{\mathop{{\rm lim\ sup}}}

  %%Fractions

\def \frac#1#2{{ {#1\over #2}}}
\def \fracd#1#2{{\displaystyle {#1\over #2}}}

\def \qed {\hfill$\square$\par}

%% Special

%\def\no {\noindent}
\def\=d{{\,\buildrel (d) \over =\,}}
\def\a.s.{{\buildrel a.s. \over \longrightarrow}}
\def\ra{\rightarrow}

\def\proof{\med {\it Proof. }}
\def\tf{\wt f}

\def\OI{{\mathop {{\rm OI}}}}
\def\UC{{\mathop {{\rm UC}}}}

\font\tf=cmbx10 scaled\magstephalf
\font\smf=cmr8

%START
\cl { \tf Some remarks on the elliptic Harnack inequality}

\vskip 0.4 truein 

\centerline {Martin T. Barlow\footnote{$^1$}{Research partially
 supported by grants from NSERC (Canada) and CNRS (France).}}
\centerline {Department of Mathematics}
\centerline {University of British Columbia }
\centerline {Vancouver, V6T 1Z2}
\centerline {Canada} 

\vskip 0.3 truein
\cl {\it Abstract}

{\smf In this note we give three short results concerning the elliptic
Harnack inequality (EHI), in the context of random walks on
graphs. The first is that the EHI implies polynomial growth of the
number of points in balls, and the second that the EHI is equivalent
to an annulus type Harnack inequality for Green's functions. The third
result uses the lamplighter group to
give a counterexample concerning the relation of coupling with the EHI.}

\bigb {\it Keywords:} Graph, random walk, elliptic Harnack inequality,
coupling, lamplighter group.

\vskip 0.3 truein
The study of the relation between geometric properties of spaces
(in particular, manifolds), and the properties of heat kernels
and harmonic functions on these spaces has a long history, and
continues to be very active. In particular, the papers
[FS], [Gg], [SC] establish the equivalence between 
the parabolic Harnack inequality (PHI),  
Gaussian type estimates on the heat kernel, 
and the condition that the underlying space satisfies volume
doubling and a family of weak Poincar\'e inequalities.
These results were translated into the graph context in [D1].

\smallskip
This characterizes PHI, and imply that it is stable under
rough isometries.  No similar characterization is known for the weaker
elliptic Harnack inequality (EHI). It is trivial that PHI implies EHI,
while in [BB1] and [D2] examples are given of spaces which satisfy EHI
but fail to satisfy the PHI.  The example in [D2] shows in addition
that EHI does not imply volume doubling. See [HSC] for some further
results on the difference between EHI and PHI.

In this note we give three short results concerning the EHI. For
simplicity we give these in the easiest context, that of random
walks on graphs. The first is that EHI implies polynomial
growth of the number of points in balls, and the second that EHI
is equivalent to an annulus type Harnack inequality for Green's
functions. The final result uses the lamplighter group to show
that a coupling type condition for random walks does not imply
EHI. The paper is concluded with an open problem.

\medskip 
We begin by recalling the definition of a weighted graph.
Let $\Gam=(G,E)$ be an locally finite connected graph.  If $x$ and $y$
are connected by an edge $\{x,y\}$ we write $x \sim y$.  We call
$\nu=(\nu_{xy})$, $x,y\in G$ a {\sl conductance matrix} if
$\nu_{xy}\ge 0$, $\nu_{xy}=\nu_{yx}$ for all $x,y \in G$ and in
addition $\nu$ is linked to the graph structure by the condition that
$\nu_{xy}>0$ if and only if $x\sim y$. The triple  $(G,E, \nu)$ is called
a {\sl weighted graph}.  We call the {\sl natural weight} on $\Gam$
the weights given by taking $\nu$ to be the adjacency matrix of
$\Gam$; that is $\nu_{xy} = 1$ if and only if $x\sim y$.
Let $\mu(x) =\sum_y \nu_{xy}$;  we extend $\mu$ to a measure on $G$.
%and $\nu$  to a measure on $E$. 

Given a weighted graph $(\Gam,\nu)$ we define the simple random
walk $X$ on $G$ to be the Markov chain 
with transition probabilities given by
$$ p_{xy} = \bP^\cdot( X_{n+1}=y|X_n=x) = \fracd{\nu_{xy}}{\mu(x)}, \q
  x,y \in G, \,\, n\ge 0. \eqno(1) $$
Let $d(x,y)$ be the graph distance on $\Gam$. For $x\in G$,
$r\in (0,\infty)$, let
$$ B(x,r)=\{y: d(x,y) \le r\}.$$
Note that any pair of points $x$ and $y$ will be connected
by a geodesic path of length $n=d(x,y)$,  but that this path will not in
general be unique. We select a family  $\gamma(x,y)$, $x,y\in G$ of
geodesic paths. 

\med {\bf Definition.} 1. Let $A\subset G$. We write 
$\partial A = \{ y\in A^c: d(x,y)=1$ for some $x \in A\}$ for the exterior
boundary of $A$, and set $\ol A=A \cup \partial A$, and write
$\oB(x,r)=B(x,r)\cup\pd B(x,r)$.
\nl 2. Define the Laplacian on $(\Gam,\nu)$ by
$$ \Delta f(x) = \frac{1}{\mu_x} \sum_y \nu_{xy} (f(y)-f(x)). $$
\nl 3. A function $h$ is {\sl harmonic on} $A\subset G$ if
$h: \ol A\to \bR$ and $\Delta h(x) = 0$, $x\in A$.

\smallskip
We now introduce two conditions that a weighted graph
$(\Gam,\nu)$ may or may not satisfy.

\med {\bf Definition.} 1.
$(\Gam,\nu)$ has {\sl controlled weights}
if there exists $p_0>0$ such that for all $x\in G$ 
$$ \frac{\nu_{xy}}{\mu(x)} \ge p_0. \eqno(2) $$
(This was called the $p_0$-condition in [GT2].)

\med 2.  $(\Gam,\nu)$ satisfies {\sl an elliptic Harnack inequality} (EHI) if
there exists $C_1>0$ such that, for any $x\in G$, $R\ge 1$, and
non-negative $h:G \ra \bR$ harmonic in $B(x,2R)$,
$$ \sup_{B(x,R)} h \le C_1  \inf_{B(x,R)} h. \eqno(3) $$
We have taken balls $B(x,R)\subset B(x,2R)$ just for simplicity. 
If  $(\Gam,\nu)$ has controlled weights, 
$K>1$ and (3) holds whenever $h\ge 0$ is harmonic in $B(x,KR)$,
then an easy chaining argument gives (EHI) (for a different constant $C_1$).

\ms Note that if $(\Gam,\nu)$ has controlled weights
then it satisfies a local Harnack inequality: if $h\ge 0$ is 
harmonic on $A\subset G$ then
$$ h(x) \ge p_0 h(y), \q \hbox{ if } x\in A, \q y \sim x. $$
It follows immediately that any finite graph satisfies (EHI), but
it  may still be of interest to ask how good the constant $C_1$ can be.

\med {\bf Remark.} 
A hypothesis such as controlled weights is needed to connect the
graph structure with the Laplacian. One might hope that (EHI) would
imply controlled weights, but this is not the case.
Consider the graph $\Gam$ with $G=\bZ \times \{0,1,2\}$, and
edges
$$ E = \big\{ \{(n,0),(n+1,0)\}, \{(n,i),(n,j)\}, \, n\in \bZ,
i, j=0, 1,2,\,  i \neq j \big\}. $$
Let  $\nu_{(n,1),(n,2)}=2^{-|n|}$, $n \in \bZ$, and all other 
edges have weight 1. Then $(\Gam,\nu)$ does not have controlled weights
but it is easy to verify that (EHI) holds.

\ms The first main result of this paper is that (EHI) implies a 
bound on the size of the balls $B(x,R)$. Write $|A|$ for the number
of elements in the set $A$. 

\proclaim Theorem 1. Let $(G,E,a)$ be a weighted graph with controlled
weights which satisfies (EHI) with a constant $C_1$. Then
there exist constants $C$, $\th$, depending only on $C_1$ such that
$$ |B(x_0,R)| \le C R^{1+\th}, \q x_0 \in G, R\ge 1. \eqno(4) $$

\proof Let $B=B(x_0,R)$, $\tau=\min\{n: X_n \in \pd B\}$.
Fix $z\in \pd B$, and let
$$ h(x) = h_z(x) = \bP^x( X_\tau =z). $$

Let $y_i$ be a sequence of points on the geodesic path  $\gamma(z,x_0)$
with $d(y_j,z)=3^j$, for $0\le j \le N$, where 
$N$ is chosen so that $3^N\le R < 3^{N+1}$.
Choose also $z_j$ on this geodesic with $d(z,z_j)=2 . 3^j$, for 
$0\le j\le N$, except that if $2.3^N>R$ then we take $z_N=x_0$.
Applying (EHI) to $h$ in $B(z_j, 3^j) \subset B(z_j, 2. 3^j)$
we have, if $0\le j <N$
$$ h(y_{j+1}) \ge C_1^{-1} h(y_j), $$
while if $j=N$ we obtain $h(x_0)\ge C_1^{-1} h(y_N)$.

The local HI implies that $h(y_0)\ge p_0 h(z)=p_0$, and therefore
$$ h_z(x_0) \ge \frac{p_0}{C_1} C_1^{-N} 
\ge  \frac{p_0}{C_1} R^{-\th}, $$
where $\th = \log C_1/\log 3$. Since 
$$ 1 = \sum_z h_z(x_0) \ge |\pd B| c_1 R^{-\th}, $$
we obtain  $|\pd B(x_0,R)| \le c R^\th$, and  summing over $R$
gives (4). \qed

\med {\bf Remarks.} 
1. This proof controls the number of vertices in 
$B(x_0,R)$ rather than 
 \break $V(x_0,R)=\mu(B(x_0,R))$;  I do not
know if one has similar control of $V(x_0,R)$.
\nl 2. Note that if $|B(x_0,R)|\ge R^{1+\delta}$
then this proof implies that $C_1\ge  3^{\delta}$.

\bigb We now recall the definition of the Green's function
on $(\Gam,\nu)$. Let $D\subset G$, $\tau_D=\min\{n\ge 0: X_n \not\in D\}$,
and $p_n^D(x,y)=\bP^x(X_n=y, n< \tau_D)/\mu(y)$.
Then we set 
$$ g_D(x,y) =\sum_{n=0}^\infty p_n^D(x,y).$$
Sufficient conditions for $g_D$ to be finite are that $D$ is finite,
or that $(\Gam,\nu)$ is transient.
$g_D$ is symmetric in $x$, $y$, is zero if either 
variable is outside $D$, and $g_D(x_0, \cdot)$ is harmonic in 
$D-\{x_0\}$.

In [GT2] the following condition (HG) is introduced.

\med {\bf Definition.} $(\Gam,\nu)$ satisfies  (HG) if there exists a
constant $C_2\ge 1$ such that, 
for any finite set $D\subset G$, if $R\ge 1$ and 
$B(x_0,2R) \subset D$ then
$$ \max_{y\in B(x_0,R)^c} g_D(x_0,y) \le C_2 
 \min_{y\in B(x_0,R)} g_D(x_0,y). \eqno(HG)$$

\small In [GT1] and [GT2] it is proved that (HG) implies (EHI),
and that (EHI) and another (geometric) condition,
denoted (BC), implies (HG). The point of the next result is that no
additional geometric condition is needed.

\proclaim Theorem 2. Suppose $(\Gam,\nu)$ has controlled
weights and satisfies (EHI). 
Then there exists a constant $C_2$ such that
if $x_0\in G$, $R\ge 1$, $d(x_0,x)=d(x_0,y)=R$, and $B(x_0,2R)\subset D$ 
then
$$ C_3^{-1} g_D(x_0,y) \le g_D(x_0,x) \le C_3 g_D(x_0,y). \eqno(5)$$
In particular  the conditions (EHI) and (HG) are equivalent.

\proof Note first that by symmetry it is enough to prove the right hand
inequality. Using the local Harnack inequality we can strengthen 
(EHI) to give, if $d(x_0,x)=r$, that
$$ \max_{B(x,r/2)} g_D(x_0,.) \le C_4 \min_{B(x,r/2)} g_D(x_0,.). $$
We first assume  that $R\ge 12$ and is divisible by 12. 
Let $x'$, $y'$ be the midpoints of $\gamma(x_0,x)$, and 
$\gamma(x_0,y)$. Thus $d(x_0,x')=d(x_0,y')=R/2$. Clearly we have
$d(x',y)\ge R/2$ and $d(x,y')\ge R/2$.

We now consider two cases.

\small Case 1. $d(x',y')\le R/3$. Let $z$ be as close as possible
to the midpoint of $\gamma(x',y')$.
Then $d(z,x')\le 1+R/6 \le R/4$. 
So applying (EHI) to
$g_D(x_0,\cdot)$ in $B(x',R/4)\subset B(x',R/2)$, we deduce that
$$  C_4^{-1} g_D(x_0,x')\le g_D(x_0,z) \le C_4 g_D(x_0,x') . \eqno(6)$$
Now apply (EHI) to 
$g_D(x_0,\cdot)$ in $B(x,R/2)\subset B(x,R)$, to deduce that
$$ C_4^{-1} g_D(x_0,x) \le g_D(x_0,x') \le  C_4 g_D(x_0,x). $$
Combining these inequalities we deduce that
$$  C_4^{-2} g_D(x_0,x)\le g_D(x_0,z) \le C_4^2 g_D(x_0,x),$$
and this, with a similar inequality for $g_D(x_0,y)$, 
proves (5).

\small Case 2. $d(x',y')> R/3$.
Apply (EHI) to $g_D(y,\cdot)$ in $B(x_0,R/2)\subset B(x_0,R)$, to 
deduce that
%$$ C_4^{-1} g_D(y,x') \le g_D(y,y') \le  C_4 g_D(y,x'). \eqno(7) $$
$$ C_4^{-1} g_D(y,x') \le g_D(y,x_0) \le  C_4 g_D(y,x'). \eqno(7) $$
Now look at $g_D(x',\cdot)$. 
If $z'$ is on $\gamma(y',y)$ with $d(y',z')=s\in [0, R/2]$ then
as $d(x',y')>R/3$ and $d(x',y)\ge R/2$ we have
$d(x',z') \ge \max(R/3-s, s)$. Hence we deduce $d(x',z')\ge R/6$.
So applying (EHI) repeatedly to $g_D(x',\cdot)$ for a chain of balls 
$B(z',R/12)\subset B(z',R/6)$ we deduce that
$$  C_4^{-6} g_D(x',y')\le g_D(x',y)\le C_4^6 g_D(x',y'). \eqno(8) $$
So, we obtain from (7) and (8), 
$$ g_D(y,x_0) \le  C_4 g_D(y,x') \le C_4^7 g_D(x',y'), \q
 g_D(x',y')\le C_4^6 g_D(y,x') \le C_4^7 g_D(y,x_0). $$
We have similar inequalities relating $g_D(x,x_0)$ and $g_D(x',y')$,
which proves (5).

If $R\le 11$ then (5) follows using the local Harnack inequality.
If $R\ge 12$ is not divisible by 12, let $R'=12k$ with
$R-11\le R'<R$, and choose points $x_1$ and $y_1$
a distance $R'$ from $x_0$ with $d(x,x_1)\le 11$,
$d(y,y')\le 11$. Then we have (5) for $x_1$ and $y_1$,
and (5) for $x$ and $y$ again follows using the local Harnack
inequality. 

Finally, to deduce (HG) note that by the maximum principle
the maximum and minimum in (HG) will be attained at 
a point $y$ with $d(x_0,y)=R$  and $z$ with $d(x_0,z)=R-1$.
Hence using the local Harnack inequality (HG)
follows from (5). \qed

\bigskip The final question we consider is the relation
between (EHI) and coupling for random walks on $(\Gam,\nu)$.
For further details on coupling see [Lv].
It is sometimes easier to couple continuous time random walks (CTSRW), 
rather than discrete time ones, so we recall
that the continuous time simple random walk $Y_t$ on $(\Gam,\nu)$ is the
Markov process with generator $ \Delta f(x)$.
$Y$ waits for an exponential time with mean 1 at each vertex $x$, and
then moves to each $y\sim x$ with the same probability as $X$, that
is $p_{xy}=\nu_{xy}/\mu(x)$.

\med {\bf Definition.} 1. Let $K>1$. 
$(\Gam,\nu)$ satisfies {\sl uniform (co-adapted) coupling} with constant
$K$, if for all $x_0\in G$, $R\ge 1$, and $y_1,y_2 \in B(x_0,R)$ 
there exist CTSRW walks  $Y^i_t$ with $Y^i_0=y_i$ satisfying
the following condition. Let
$$ \eqalign{
 \tau_C &=\min\{ t\ge 0: Y^1_t=Y^2_t \}, \cr
 \tau_E &=\min\{ t\ge 0: Y^1_t\notin B(x_0,KR) \hbox{ or }  Y^2_t\notin B(x_0,KR)
  \}. \cr }$$
Then there exists $p_1>0$, independent of $x_0$ and $R$, such that
$\bP(\tau_C< \tau_E) >p_1$. We denote this condition $\UC(K)$.
\nl 2.  Write  $\Osc(h,A)=\max_A h -\min_A h$.
$(\Gam,\nu)$ satisfies an {\sl oscillation inequality} with constant
$K>1$ (denoted $(\OI(K)$), if there exists $\rho<1$ such that, 
if $h$ is harmonic in $B(x,KR)$ then
$$ \Osc(h,B(x,R)) \le \rho \Osc(h,\oB(x,KR)) . $$

\med {\bf Remarks.} 1. In the above we require that the $Y^i$ are 
defined on the same filtered probability space.
\nl 2. Some of the consequences of EHI, such as 
the Liouville property for bounded  harmonic functions on $\Gam$,
also follow from $\OI(K)$.
\nl 3. There are several different types of coupling -- see 
[CG] for a weaker type.

\ms It is easy to see that uniform coupling implies $\OI(K)$.
For, if $h$ is harmonic on $B(x_0,KR)$ 
and we write $\Osc(h,A)=\max_a h -\min_A h$ then, writing $B=B(x_0,R)$,
$B^*=\ol B(x_0,KR)$, $T=\tau_C \wedge \tau_E$, and choosing $x_i$
suitably, 
$$ \eqalignno{
 \Osc(h,B)&= h(x_2)-h(x_1)  \cr
        &= \bE (h(Y^2_T)-h(Y^1_T)) \cr
  &\le \bP(\tau_E < \tau_C)  \Osc(h,B^*) 
\le (1-p_1) \Osc(h,B^*). &(9)\cr}$$ 
\ms
It is also well known that EHI implies $\OI(K)$, and that 
$\OI(K)$ plus a suitable lower bound on the hitting probabilities of small
balls implies EHI. (See for example [FS] or [BB1]).
The following example on the lamplighter group shows that this extra 
condition is needed: uniform coupling alone is not enough to imply
EHI. It is also shows that (unlike EHI) the size of the larger
ball $B(x,KR)$ plays an important role in the conditions $\UC(K)$
and $\OI(K)$.

For details of the lamplighter group see for example [LPP] or [W]. Let
$$ \sA =\{ \xi\in  \{0,1\}^{\bZ}: \sum_k \xi_k < \infty \}, $$
and $G= \bZ \times \sA$.  We denote points $x\in G$ by
$x=(n,\xi)$, and will write $A(x)=A(n,\xi)=\{i: \xi_i=1\}$.
For $\xi \in  \sA$, $k=0,1$ let 
$T_{n,k}(\xi)\in  \{0,1\}^{\bZ}$ be the sequence $\eta$ with
$\eta_i=\xi_i$, $i\neq n$, $\eta_n=k$. 
We define edges on $G$ by taking the four neighbours of $(n,\xi)$
to be the points $(n \pm 1, T_{n,k}(\xi))$, $k=0,1$.
One thinks of the $\xi_i$ as lamps (off or on), and $n$ as the location
of the lamplighter. Each move on the graph, the lamplighter may switch
the status of the current lamp, and then moves to a neighbouring point.
We write $Y_t=(U_t, \Xi(t))$ for the CTSRW on $\Gam$. Note that $U$ is a
CTSRW on $\bZ$ and that after $U$ leaves any site $n$ the state
of the lamp $\Xi_n$ has been randomized. 
(This is the `switch then walk' random walk on $G$; 
there are several other natural random walks, such as
`walk then switch' or `switch, walk, switch' -- see [W] for more details).

\proclaim Theorem 3. 
(a) The lamplighter group $G$ satisfies uniform coupling for $K>4$.
\nl (b) $G$ fails to satisfy $\OI(K)$ for $K<3$, and in particular
does not satisfy EHI.

\proof (a) By the symmetry of the group we can take $x_0=(0, \eta^0)$, 
where $\eta^0_i = 0$ for all $i$. 
If $x=(x',\xi)\in G$, with $-a=\min\{n: \xi_n=1\}$, $b=\max\{n: \xi_n=1\}$,
and $x'\ge 0$ then it is easy to check that
$d(x_0,x) = 2a + b + |b-x'|$.

Let $x_i=(n_i,\xi^i) \in B(x_0,R)$, $i=1,2$.
Then we clearly have $|n_i|\le R$ and 
$\{n: \xi^i_n =1 \} \subset [-R, R]$.
Let $0< \eps < (K-4)/8$. 
We now define $Y^i_t$ with initial states $x_i$ as follows. 
First, we try to couple the CTRSW $U^1$ and $U^2$ using
reflection coupling -- see [LR]. 
With probability $p_1=p_1(\eps)>0$, independent of $R$, 
this coupling  succeeds at a time $T_1$
before either $U^1$ or $U^2$ has left $I(\eps)=[-R(1+\eps),R(1+\eps)]$. 
Denote this event $F_1$. 

On $F_1$ we have $\{n: \Xi^i_n(T_1) =1 \} \subset I(\eps)$.
If $F_1$ occurs, we now continue by taking $U_t=U^1_t=U^2_t$
for $t\ge T_1$, and by making each lamp randomization the same
for the two processes. Hence $\Xi^i_n(t)$ agree for every site $n$
visited by $U$ between $T_1$ and $t$.

Let $F_2$ be the event that $U$ hits both $\pm R(1+\eps)$
before leaving $I(2\eps)$; then $\bP(F_2|F_1) >p_2>0$. 
Let $T_2$ be the exit time from $I(2\eps)$.
Then 
$\bP(F_1\cap F_2)\ge p_1p_2$ and on $F_1\cap F_2$
we have $Y^1_{T_2}=Y^2_{T_2}$, and 
 $\{n: \Xi^i_n(t) =1 \} \subset I(2\eps)$, $U^i_t \in  I(2\eps)$,
for each $t\in [0, T_2]$. Thus on $F_1\cap F_2$
$d(x_0,Y^i_{t}) \le 4(1+2\eps)R < KR$ for $t\in [0, T_2]$,
so that $T_2 < \tau_E$.
This proves that $\Gam$ satisfies uniform coupling with 
constant $K$.

\med (b) We now prove that OI with $K<3$ fails for $G$.
For this calculation it is easier to use the discrete
time random walk $X_n=(V_n,\Theta(n))$; of course $X$
and $Y$ have the same harmonic functions.
Note that if $T$ is a stopping time with respect to
$V$, then at time $T$ the states of all the lamps
$\Theta_n(T)$, $n\in \{V_k, 0\le k\le T-1\}$ will have been
randomized.

Let $\delta= (3-K)/3$, $\lam=1-\delta$.
%For $y=(u,\xi)\in G$ let $I(y)=\{n: \xi_n=1\}$.
Let $\tau=\min\{n: d(x_0,X_n)>KR\}$,
set $G =\{ A(X_\tau) \cap [\lam R, R]\neq \emptyset \}$
and let 
$$ h(y) = \bP^y (G), \q y\in B(x_0,KR). $$
Let $y_1=(-R, \xi^{(1)})$, $y_2=(R, \xi^{(2)})$, where
$\xi^{(1)}_n = 1_{[-R,0]}(n)$ and 
$\xi^{(2)}_n = 1_{[0,R]}(n)$.

We begin by bounding $h(y_1)$. Let $V_0=-R$, and
$T= \min\{n: V_n \not\in (-3R,\lam R) \}$. Note that
since $A(X_T)\cap[\lam R,R]=A(X_0)\cap[\lam R,R]$, 
we have $G \subset \{ T < \tau\}$.
So 
$$ h(y_1) = \bP^{y_1}(G \cap \{T < \tau\})
\le  \bP^{y_1}(T < \tau). $$
Let $-Z=\min\{ n: \Theta_n(T)=1\}$. 
If $T<\tau$ then $V_T=\lam R$, and 
we have 
$$ d(x_0,X_T) = 2Z + V_T = 2Z + \lam R < K R, $$
so that $Z\le \lam R$. Thus at time $T$ every lamp
in $[-R, -\lam R)$ must have been switched off, and we
deduce that
$$ \bP^{y_1}(G) \le 2^{-\delta R}. $$

\ms Similar considerations lead one to expect that
$\bP^{y_2}(G^c) \le c 2^{-\delta R}$, but 
the argument is a little more involved. A cruder 
bound is easier to obtain, and enough here.
Let $V_0=R$, and note that, for any $n\ge 0$,
$$ \bP^{y_2}( A(X_n) \cap (\lam R, R]=\emptyset) = 2^{-\delta R}. $$
Let $T$ be the first hit by $V$ on $\pm 3R$: we have
$T> \tau$. So
$$ \eqalign{
 \bP^{y_2}(G^c)&\le  \bP^{y_2}(T \ge R^3) +  \bP^{y_2}(G^c;V< R^3)\cr
 &\le {c}{R^{-1}} +  \bP^{y_2}( A(X_n)\cap (\lam R,R]=\emptyset
 \hbox{ for some } 0\le n< R^3) 
 \le {c}{R^{-1}} + R^3 e^{-\delta R} . }$$
So $\Osc(h, \oB(x_0,R)) \ge h(y_2)-h(y_1) = 1 -o(R)$,
and hence $\OI(K)$ fails for $G$. 

Finally, since  EHI implies OI for any $K>1$,   
EHI also fails for $G$.
\qed

\med {\bf Remarks.} 1. [W] identifies the cone of positive harmonic 
functions associated with the `walk then switch' and the
`switch, walk, switch' walks on $G$: in both cases there are non-constant
positive harmonic functions.
\nl 2. The key feature used by the argument in Theorem 3 
is that after $U$ has moved across an interval $I\subset \bZ$
all the lamps in $I$ will have been randomized. So this proof also
works for the `walk then switch' and `switch, walk, switch' walks.
But it does not apply to the `walk or switch' walk.
\nl 3. It seems likely that there is a gap between the best constant
$K$ for $\UC(K)$ and $\OI(K)$.

\ms We conclude this paper by mentioning an open problem. 
If $D \subset G$ we define for $f:D \to \bR$ the Dirichlet form
on $L^2(D)$
$$ \sE_D(f,f)= \sum_{x\in D} \sum_{y\in D}\nu_{xy}(f(x)-f(y))^2. $$
For disjoint $A, B \subset D$ define the conductance 
$$  C_D(A,B) =\inf \{ \sE_D(f,f): f= 1 \hbox { on $A$ and $f=0$ on $B$}\}. $$
Given $x_0\in G$, $R\ge 10$ let 
$$ \sD(x_0,R)=\{(x,y): d(x,x_0)\le R/2, d(y,x_0)\le R/2, d(x,y)\ge R/3\}. $$

We introduce the following condition.

\med {\bf Definition.} $(\Gam,\nu)$ satisfies the 
{\sl dumbbell condition} (DB) if there exists a constant $C_3<\infty$ such that
for all $x_0\in G$, $R\ge 10$, writing $D=B(x_0,R)$, $\sD= \sD(x_0,R)$,
$$ \max_{(x,y)\in \sD} C_D(B(x,R/10),B(y,R/10))
 \le C_3  \min_{(x,y)\in \sD} C_D(B(x,R/10),B(y,R/10)). \eqno(DB)$$
Note that (DB) is stable under bounded perturbation of the
weights $\nu_{xy}$: if  $(\Gam,\nu)$ satisfies (DB) and
$\nu'_{xy}$ satisfy, for some constant $c_1$,
$$ c_1^{-1} \nu_{xy} \le \nu'_{xy} \le  c_1 \nu_{xy}, \q x,y \in G, $$
then   $(\Gam,\nu')$ also satisfies (DB).

\med {\bf Problem.} Suppose  $(\Gam,\nu)$ has controlled weights.
Is (DB) equivalent to EHI?

\med {\bf Acknowledgment.} I am grateful for Terry Lyons for 
some stimulating conversations related to this paper.

\bigskip \cl { \bf References}

\nr [BB1] M.T. Barlow and R.F. Bass, Brownian motion and harmonic analysis on
Sierpinski carpets. {\it Canad. J. Math.} {\bf 51} (1999) 673--744.

\nr [CG] M. Cranston, A. Greven. Coupling and harmonic functions in the case
of continuous time markov processes. {\it Stoch. Proc. Appl.} {\bf 60}
(1995), 261--286.
   
\nr [D1] T. Delmotte. Parabolic Harnack inequality and estimates
of Markov chains on graphs.  {\it Rev. Math. Iberoamericana} {\bf 15}
(1999), 181--232.

\nr [D2] T. Delmotte. Graphs between the elliptic and parabolic Harnack
inequalities.  {\it Potential Anal.} {\bf 16} (2002), 151--168. 

\nr [FS] E.B. Fabes and D.W. Stroock. A new proof of Moser's parabolic
Harnack inequality via the old ideas of Nash.
{\it Arch. Mech. Rat. Anal. \bf 96} (1986) 327--338.

\nr [Gg] A.A. Grigor'yan. The heat equation on noncompact Riemannian
manifolds. {\it Math. USSR Sbornik \bf 72} (1992) 47--77.

\nr  [GT1] A. Grigor'yan, A. Telcs. Sub-Gaussian estimates of heat kernels
on infinite graphs. {\it Duke Math. J.} {\bf 109} (2001) 452--510.

\nr [GT2]  A. Grigor'yan, A. Telcs. Harnack inequalities and sub-Gaussian 
estimates for random walks. {\it Math. Annalen} {\bf 324} (2002), 521--556. 

\nr [HSC] W. Hebisch, L. Saloff-Coste. On the relation between
elliptic and parabolic Harnack inequalities. 
{\it Ann. Inst. Fourier (Grenoble)} {\bf  51}  (2001), 1437--1481. 

\nr [Lv] T. Lindvall, {\it Lectures on the coupling method}.
Wiley, New York, 1992.

\nr [LR] T. Lindvall and L.C.G. Rogers. Coupling of multi-dimensional diffusions by 
reflection. {\it Ann. Prob. \bf 14} (1986) 860--872.

\nr [LPP] R. Lyons, R. Pemantle, Y. Peres.  Random walks on the lamplighter group. 
{\it Ann. Probab.} {\bf 24} (1996), no. 4, 1993--2006.

%\nr [M1]  J. Moser. On Harnack's inequality for elliptic differential equations.
%{\it Comm. Pure Appl. Math.} {\bf 14} (1961) 577--591.

%\nr [M2] J. Moser. On Harnack's inequality for parabolic differential equations.
%{\it Comm. Pure Appl. Math. \bf 17} (1964) 101--134. 

\nr [SC] L. Saloff-Coste. A note on Poincar\'e, Sobolev, and Harnack
inequalities. {\it Duke Math. J. } {\bf 65} (1992), 
{\it Inter. Math. Res. Notices} {\bf 2}  (1992) 27--38.

\nr [W] W. Woess. Lamplighters, Diestel-Leader graphs, random walks,
and harmonic functions. Preprint 2003.

\med \ver \q \vdate

\bye